\documentclass{article}
\usepackage[utf8]{inputenc}
\usepackage{graphicx}
\usepackage{caption}
\usepackage{subcaption}

\title{Discussion of ``A note on universal inference''\\ by Timmy Tse and Anthony Davison}
\author{Mathias Drton, Hongjian Shi, David Strieder}
\date{December 2022}

\begin{document}

\maketitle

We congratulate the authors for their insightful contribution to universal inference, as proposed in \cite{MR4242731}.  The split likelihood ratio test (SLRT) provided by the universal inference framework is valid in finite samples and requires essentially no regularity conditions.  However, this amazing universality also comes at a price as the SLRT may be highly conservative.

In their article, Tse and Davison \cite{tse2022note} 
explore how nuisance parameters impact the conservativeness of the SLRT. They study the loss of power in the basic normal location model, taking up the analysis in \cite{dunn2022gaussian}.  
Tse and Davison also propose a modification to reduce the conservatism of the SLRT.
Universal inference is justified by an application of Markov's inequality, and the modification is based on the clever observation that, in the normal location model, the expectation to be considered may be more tightly controlled by exploiting the independence of two terms that emerge in the presence of nuisance parameters.
The impact of using this trick for irregular testing problems is illustrated in an example based on a Gaussian mixture model.  

An independent recent work of our group studies universal inference from a related but also different angle \cite{strieder2022choice}.  We would like to take the opportunity of this discussion to point out similarities but also differences between this work and the paper of Tse and Davison.

The main results in Strieder and Drton \cite{strieder2022choice} provide the large-sample theory for the SLRT under both the null hypothesis and local alternatives assuming the classical regularity conditions that lead to Wilks' Theorem for the standard likelihood ratio test (see Theorem~3.1 and related Corollaries/Remarks in \cite{strieder2022choice}).
The setup allows for general composite null hypotheses.
The limiting distributions of the SLRT are then used to illustrate how (loss of) power of the test is impacted by dimensionality of the hypotheses.  These limiting distributions coincide with the distributions that are considered in the paper of Tse and Davison (e.g., Equations (8), (13)--(14)).

The two articles differ in the way the (limiting) distribution of the SLRT is used in order to improve power.
While Tse and Davison focus on using the independence of the two parts stemming from the parameter of interest and the nuisance parameter, Strieder and Drton concentrate on highlighting the impact of the data splitting. 
Tse and Davison utilize insights from the regular setting under the perspective of nuisance parameters to motivate the use of a smaller critical value.  Such a smaller critical value  obviously increases power and continues to control the size in the regular case.  However, with the modified critical value one no longer has a universal guarantee to control the size in general irregular problems.  This raises the question:  How far should we go in terms of lowering the critical value once we start setting it based on asymptotic distributions derived for regular testing problems?  While this would be easy to do, we certainly do not want to go all the way to the relevant quantile of the regular limiting distribution.  But then which principles should we invoke to justify a particular lowered critical value when applying it to irregular problems?

As mentioned, the work of Strieder and Drton pursues a different goal and proposes a different modification.  Instead of changing the critical value, this modification  employs the regular limiting distribution merely to tune the data splitting ratio (at a target power level).  Tuning the splitting ratio has a more limited impact on power, but importantly an SLRT with tuned splitting remains universally valid in finite samples---no matter whether the problem is regular or not.

\begin{figure}
\centering
\begin{subfigure}{.49\textwidth}
  \centering
  \includegraphics[width=.95\linewidth]{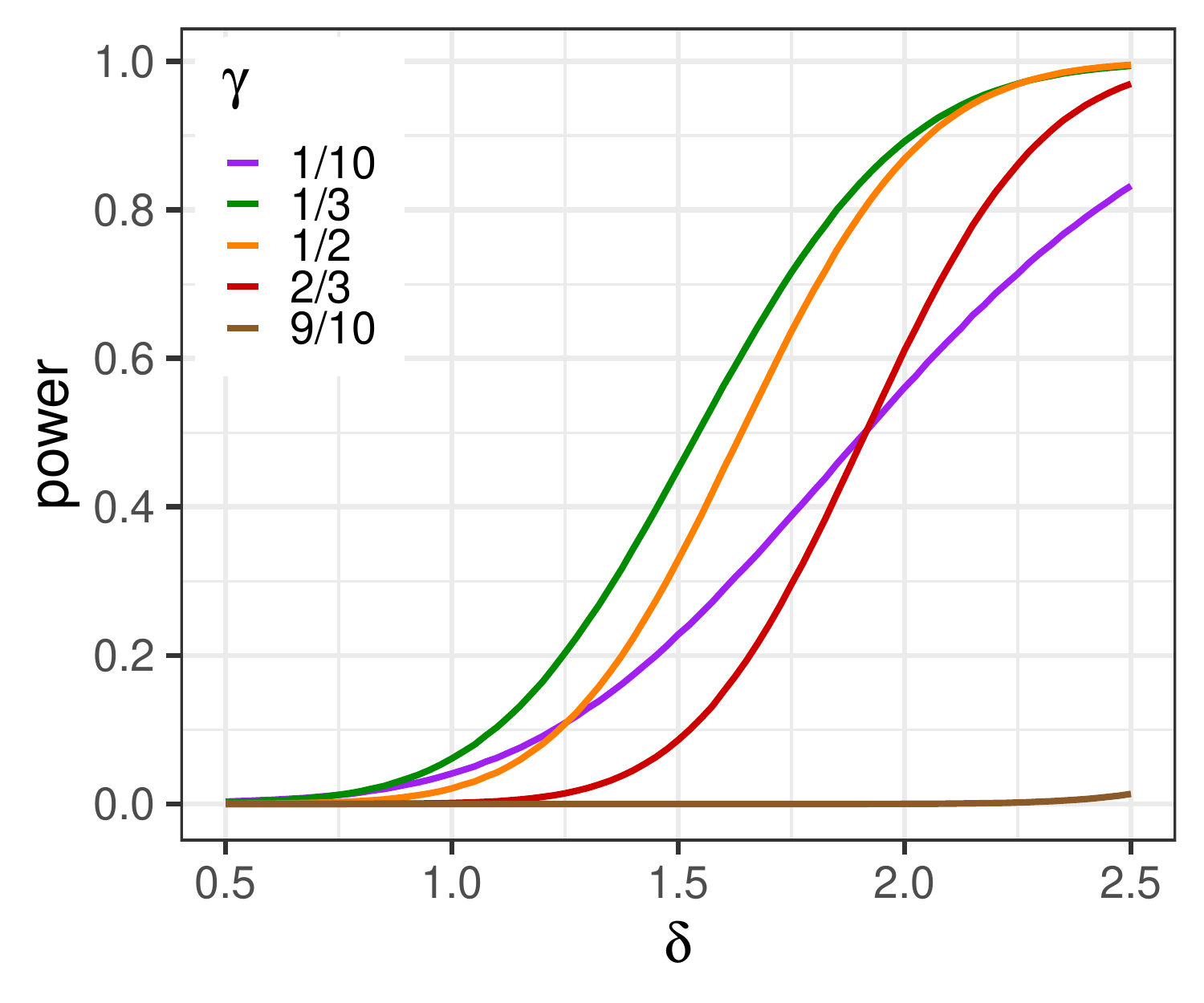}
  \caption{Power as a function of $\delta$}
\end{subfigure}
\begin{subfigure}{.49\textwidth}
  \centering
  \includegraphics[width=.95\linewidth]{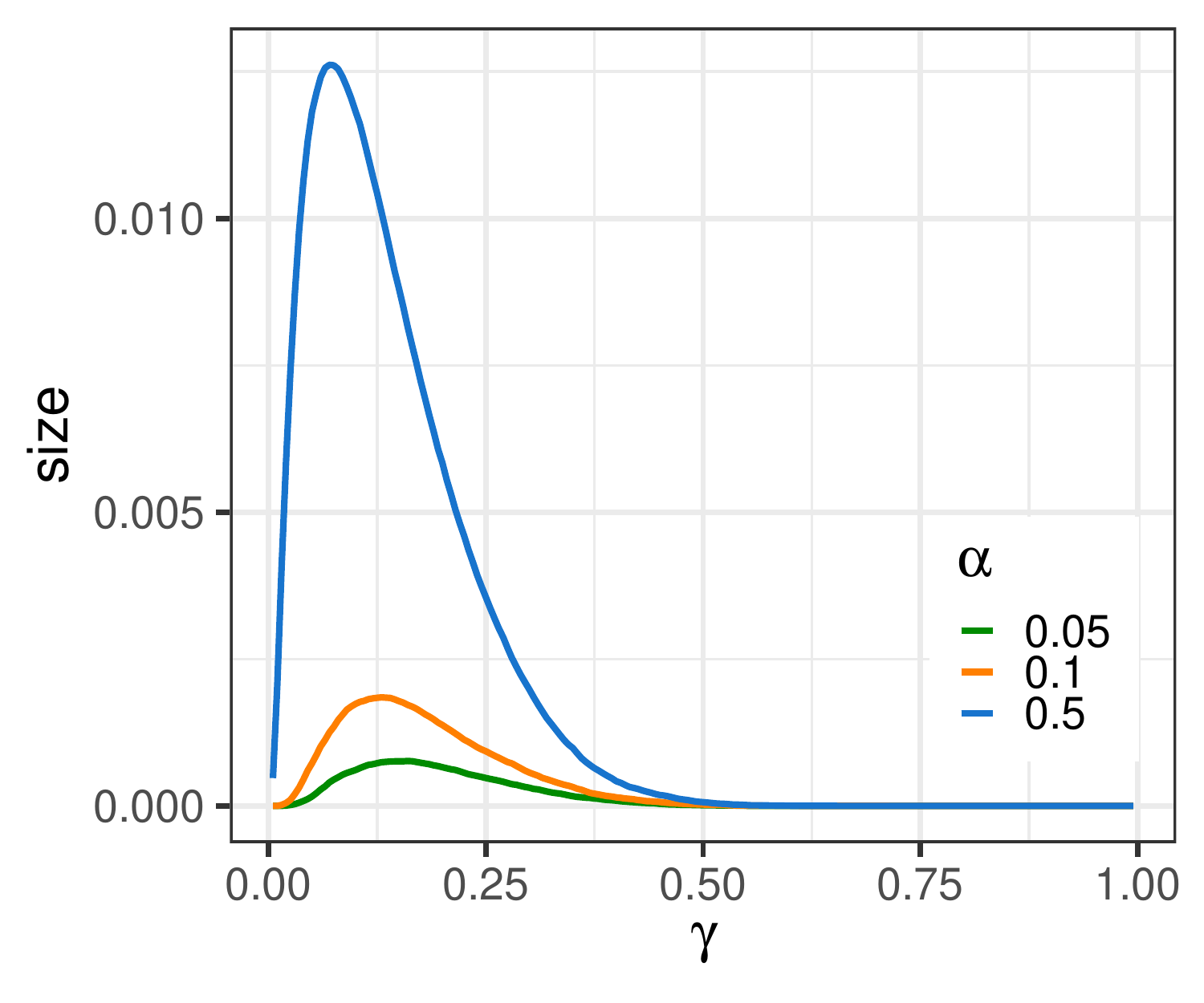}
  \caption{True size as a function of $\gamma$}
\end{subfigure}
\caption{Power of the SLRT for testing a $45$-dim.~hypothesis in $50$-dim.~parameter space with departure $(\delta,\dots,\delta)$, split ratio $\gamma$ and nominal size $\alpha$.}
\label{figure}
\end{figure}

Tse and Davison also investigate the power of the SLRT across various splitting ratios with the conclusion that a ratio of $\gamma=2/3$ performs best. We emphasize that this is a feature of the considered low-dimensional setting. Contrasting Figure 2 from the article of Tse and Davison, Figure~\ref{figure} here displays the power and size of the SLRT when testing a $45$-dimensional hypothesis in a $50$-dimensional space and shows that in a higher-dimensional setting a splitting ratio below $0.5$ performs best. For a detailed analysis of the optimal choice of the splitting ratio based on the dimensionality of tested null and alternative hypotheses, we refer to \cite{strieder2022choice}. 

\bibliographystyle{amsplain}
\bibliography{discussion_universal_inference.bib}

\end{document}